\documentclass{amsproc}
\usepackage{amssymb}
\usepackage{amsmath}
\usepackage{amsthm}
\usepackage[all]{xy}
\usepackage[ansinew]{inputenc} 
\usepackage{graphicx}

\newcommand{\tl}{\widetilde}
\newcommand{\Ol}{{\mathcal O}}
\newcommand{\f}{\varphi}
\newcommand{\ac}{{\mathcal H}}

\newcommand{\V}{{\mathcal V}}
\newcommand{\cW}{{\mathcal W}}

\newcommand{\G}{{\mathbb G}}

\newcommand{\pu}{{\mathbb P^1}}
\newcommand{\proj}{\mathbb P}
\newcommand{\quadr}{\mathbb Q}

\newcommand{\pd}{{\mathbb P^2}}

\DeclareMathOperator{\loc}{\mathrm{Locus}}
\DeclareMathOperator{\cloc}{\mathrm{ChLocus}}
\DeclareMathOperator{\cone}{{NE}}
\DeclareMathOperator{\pic}{Pic}
\DeclareMathOperator{\cycl}{N_1}

\DeclareMathOperator{\Exc}{Exc}

\newcommand{\onespan}[1]{\langle #1 \rangle}
\newcommand{\twospan}[2]{\langle #1,#2 \rangle}
\newcommand{\threespan}[3]{\langle #1,#2,#3 \rangle}
\newcommand{\morespan}[2]{\langle #1, \ldots, #2 \rangle}

\newcommand{\ms}{\medskip}

\newcommand{\om}{\textrm{Hom}}
\newcommand{\Aut}{\textrm{Aut}}

\newcommand{\ratcurves}{\textrm{Ratcurves}^n(X)}

\newcommand{\rc}[2]{#1 \xymatrix{\ar@{-->}[r] & }{#2}}

\newtheorem{theorem}{Theorem}[section]
\newtheorem{lemma}[theorem]{Lemma}
\newtheorem{proposition}[theorem]{Proposition}
\newtheorem{corollary}[theorem]{Corollary}

\theoremstyle{definition}
\newtheorem{definition}[theorem]{Definition}
\newtheorem{statement}[theorem]{}
\theoremstyle{remark}
\newtheorem{remark}[theorem]{Remark}

\begin{document}

\begin{abstract}
We study Fano manifolds of pseudoindex greater than one and dimension greater than five, 
which are blow-ups of smooth varieties along smooth centers  of dimension
equal to the pseudoindex of the manifold.\\
We obtain a classification of the possible cones of curves of these manifolds, and
we prove that there is only one such manifold without a fiber type elementary contraction.
\end{abstract}
\pagestyle{plain}

\title{Fano manifolds and blow-ups of low-dimensional subvarieties}
\author{Elena Chierici}
\author{Gianluca Occhetta}

\maketitle

\section{Introduction}
A smooth complex projective variety $X$ is called {\sf Fano} if its
anticanonical bundle $-K_X$ is ample; the {\sf index} $r_X$ of $X$ 
is the largest natural number $m$ such that $-K_X=mH$ for some (ample) 
divisor $H$ on $X$, while the {\sf pseudoindex} $i_X$ is the minimum 
anticanonical degree of rational curves on $X$.\par
By the Cone Theorem the cone $\cone(X)$ generated by the numerical classes 
of irreducible curves on a Fano manifold $X$ is polyhedral.
By the Contraction Theorem to each extremal ray of $\cone(X)$ is associated
a contraction, i.e. a proper morphism with connected fibers onto a normal variety.\\
A natural question which arises from the study of Fano manifolds 
is to investigate - and possibly classify - Fano manifolds which admit 
an extremal contraction with special features: for example, this has been done
in many cases in which the contraction is a projective bundle 
\cite{SzW1,SzW2,SzW3,SzW4,APW,NO}, a 
quadric bundle \cite{Wi} or a scroll \cite{BalW,Lan}.\par
Recently, Bonavero, Campana and Wi\'sniewski have considered the case where an extremal
contraction of $X$ is the blow-up of a smooth variety along a point, giving
a complete classification \cite{BCW}.
The case where the center of the blow-up is a curve has shown to be much more complicated. 
A complete classification in case $i_X \ge 2$ has been obtained in \cite{AOlong}, following
a more general theorem, where the classification of Fano manifolds
with a contraction which is the blow-up of a manifold along a smooth subvariety 
of dimension $\le i_X - 1$ is achieved.
As for Fano manifolds of pseudoindex $i_X=1$ which are blow-ups of smooth varieties
along a smooth curve, some special cases  
have been dealt with in the PhD thesis of Tsukioka \cite{Tsuk} 
(partially published in \cite{Tsuk2}). \par
Considering the case when the dimension of the center of the blow-up is $i_X \ge 2$, 
the lowest possible dimension of the manifold is five; 
the cones of curves of such varieties are among those listed in the in \cite{CO},
where the cone of curves of Fano manifolds of dimension five and pseudoindex greater than one
were classified. Under the stronger assumption that $r_X \ge 2$ the complete list of
Fano fivefolds which are blow-ups of smooth varieties along smooth surfaces has been 
given in \cite{CO2}.\par
\medskip
In this paper we propose a generalization of both the results in \cite{AOlong} and in
\cite{CO2}, considering Fano manifolds of dimension greater than five 
with a contraction which is the blow-up of a manifold along a smooth subvariety 
of dimension $i_X \ge 2$.\\
We will first give a classification of the possible cones of curves of these varieties:

\begin{theorem}\label{main} 
Let $X$ be a Fano manifold of pseudoindex $i_X \ge 2$ and dimension $n \ge 6$, 
with a contraction $\sigma\colon X \to Y$, associated to an extremal ray $R_\sigma$,
which is a smooth blow-up with center a smooth subvariety $B$ of dimension
$\dim B = i_X$.\\
Then the possible cone of curves of $X$ are listed in the following table,
where  $F$ stands for a fiber type 
contraction and $D_{n-3}$ for the blow-up of a smooth variety
along a smooth subvariety of codimension three.\par
\medskip
\begin{center}
\begin{tabular}{c|c||c|c|c|c|c}
\quad $\rho_X$ \quad & \quad $i_X$ \quad &  \quad $R_1$
\quad & \quad $R_2$ \quad & \quad $R_3$ \quad & \quad $R_4$ \quad & \\
\hline\hline
  2 & & $R_\sigma$ & $F$ & & & (a)\\
  2 &  & $R_\sigma$ & $D_{n-3}$  & &  & (b)\\
  3 & 2,3  & $R_\sigma$ & $F$ & $F$ &  & (c)\\
  3 & 2 & $R_\sigma$ & $F$ & $D_{n-3}$ &  & (d)\\
  4 &  2   & $R_\sigma$ & $F$  & $F$ & $F$ & (e)\\
\hline
\end{tabular}
\end{center}~\par
\end{theorem}

We will then prove that there is only one Fano manifold satisfying the assumption
of Theorem \ref{main} whose cone of curves is as in case (b) - or, equivalently, 
which does not admit a fiber type contraction: 

\begin{theorem}\label{grass}
Let $X$ be a Fano manifold of dimension $n \ge 6$ and pseudoindex $i_X \ge 2$,
 which is the blow-up of another Fano manifold $Y$ along a smooth subvariety $B$ of dimension $i_X$;
assume that $X$ does not admit a fiber type contraction.\\
Then $Y \simeq \G(1,4)$ and $B$ is a plane of bidegree $(0,1)$. 
\end{theorem}

We note that, in view of the classification given in Theorem \ref{main}
Generalized Mukai conjecture \cite{BCDD,ACO} holds for the Fano manifolds
we are considering.\\

Let us point out that the assumption $i_X \ge 2$ is essential for our methods, as well as
for the ones used in \cite{AOlong}, \cite{CO} and \cite{CO2}, on which they are based.\\

The proofs of Theorems \ref{main} and \ref{grass} are contained in section 5 and 6.
In section five we consider manifolds which possess a quasi-unsplit dominating family,
proving that they are as in Theorem \ref{main}, cases (a) and (c)-(e).\\
In section six we consider manifolds which do not possess a family as above,
proving first that their cone of curves is as in case (b), and then that the only
manifold is  the blow-up of $\G(1,4)$ along a plane of bidegree $(0,1)$.

\section{Background material}
\subsection{Fano-Mori contractions}

Let $X$ be a smooth Fano variety of dimension $n$ and let $K_X$ be its canonical divisor.
By Mori's {\sf Cone Theorem} the cone $\cone(X)$ of effective 1-cycles, which is contained in
the $\mathbb R$-vector space $N_1(X)$ of 1-cyles modulo numerical equivalence,
is polyhedral; a face $\tau$ of $\cone(X)$ is called an 
{\sf extremal face}  
and an extremal face of dimension one is called an {\sf extremal ray}.\\
To every extremal face $\tau$ one can associate a morphism $\f:X \to Z$ with connected fibers 
onto a normal variety; the morphism $\f$ contracts those curves whose numerical class lies in $\tau$, 
and is usually called the {\sf Fano-Mori contraction}
(or the {\sf extremal contraction}) associated to the face $\tau$.
A Cartier divisor $D$ such that $D = \f ^*A$ for an ample divisor $A$ on $Z$
is called a {\sf supporting divisor} of the map $\f$ (or of the face $\tau$).\\
An extremal ray $R$ is called {\sf numerically effective}, or of
{\sf fiber type}, if $\dim Z < \dim X$, otherwise the ray is {\sf non nef} or {\sf birational};
the terminology is due to the fact that if $R$ is non nef  there exists an irreducible divisor 
$D_R$ which is negative on curves in $R$.
We usually denote with $E = E(\f):= \{ x \in X\ |\  \dim\f^{-1} (\f(x)) > 0\}$
the {\sf exceptional locus} of $\f$; if $\f$ is of fiber type then of course $E=X$.
If the exceptional locus of a birational ray $R$ has codimension one,
the ray and the associated contraction are called {\sf divisorial}, otherwise they are called {\sf small}.

\subsection{Families of rational curves}

For this subsection our main reference is \cite{Kob}, 
with which our notation is coherent; for missing proofs and details see also \cite{ACO} and \cite{CO}.\\
Let $X$ be a normal projective variety and let $\om(\pu,X)$ be the scheme parametrizing
morphisms $f\colon \pu \to X$; we consider the open subscheme $\om_{bir}(\pu,X) \subset \om(\pu,X)$, 
corresponding to those morphisms which are birational onto their image, and its normalization
$\om^n_{bir}(\pu,X)$; the group $\Aut(\pu)$ acts on $\om^n_{bir}(\pu,X)$ and the
quotient exists.\par

\smallskip
\begin{definition}
The space $\ratcurves$ is the quotient of $\om^n_{bir}(\pu,X)$ by $\Aut(\pu)$;
we define a {\sf family of rational curves} to be an irreducible component
$V \subset \ratcurves$.\\
Given a rational curve $f\colon\pu \to X$ we will call a {\sf family of
deformations} of $f$ any irreducible component $V \subset
\ratcurves$ containing the equivalence class of $f$.
\end{definition}

\medskip
Given a family $V$ of rational curves, we have the following basic diagram,
where $p$ is a $\pu$-bundle induced by the projection 
$\om^n_{bir} (\pu, X)\times \pu \to \om^n_{bir} (\pu, X)$ and 
$i$ is the map induced by the evaluation $ev\colon\om^n_{bir} (\pu, X)\times \pu \to X$
via the action of $\Aut(\pu)$:
$$
\xymatrix{p^{-1}(V)=:U  \ar[r]^(.65){i} \ar[d]_{p} & X\\
 V & & }$$\\

We define $\loc(V)$ to be the image of $U$ in $X$;
we say that $V$ is a {\sf dominating family} if $\overline{\loc(V)}=X$.

\begin{remark}\label{avoids} If $V$ is a dominating family of rational curves, then its general
member is a free rational curve. In particular, by \cite[II.3.7]{Kob}, if $B$ is a subset
of $X$ of codimension $\ge 2$, a general curve in $V$ does not meet $B$.
\end{remark}

\begin{corollary}\label{E0} Let $\sigma\colon X \to Y$ be a smooth blow-up with center $B$ of
codimension $\ge 2$ and exceptional locus $E$, let $V$ be a dominating family of 
rational curves for $Y$ and let $V^*$ be a family of deformations of the strict transform
of a general curve in $Y$. Then $E \cdot V^*=0$. 
\end{corollary}

For every point $x \in \loc(V)$, we will denote by $V_x$ the 
subscheme of $V$ parametrizing rational curves passing through $x$.\par

\begin{definition} 
Let $V$ be a family of rational curves on $X$. We say that
\begin{itemize}
\item $V$ is {\sf unsplit} if it is proper; 
\item $V$ is {\sf locally unsplit} if every component of $V_x$ is proper for the general $x \in \loc(V)$.
\end{itemize}
\end{definition}

\begin{proposition}\label{iowifam}{\rm\cite[IV.2.6]{Kob}} 
Let $X$ be a smooth projective variety and $V$ an unsplit family of rational curves.
Then for every point $x \in \loc(V)$ we have
 \begin{itemize}
      \item[\rm(a)] $\dim X -K_X \cdot V \le \dim \loc(V)+\dim \loc(V_x) +1;$
      \item[\rm(b)] $-K_X \cdot V \le \dim \loc(V_x)+1$.
   \end{itemize}
\end{proposition}

In case $V$ is the unsplit family of deformations of an extremal
rational curve of minimal degree, Proposition \ref{iowifam} gives the {\sf fiber locus inequality}:

\begin{proposition}\cite{Io, Wicon}\label{fiberlocus} 
Let $\f$ be a Fano-Mori contraction
of $X$ and $E$ its exceptional locus$;$
let $F$ be an irreducible component of a (non trivial) fiber of $\f$. Then
$$\dim E + \dim F \geq \dim X + l -1$$
where $l =  \min \{ -K_X \cdot C\ |\  C \textrm{~is a rational curve in~} F\}.$
If $\f$ is the contraction of an extremal ray $R,$ then $l$ is called the {\sf length of the ray}.
\end{proposition}

\begin{definition}\label{CF}
We define a {\sf Chow family of rational curves} $\V$ to be an irreducible 
component of  $\textrm{Chow}(X)$ parametrizing rational and connected 1-cycles.\\
If $V$ is a family of rational curves$,$ the closure of the image of
$V$ in $\textrm{Chow}(X)$ is called the {\sf Chow family associated to} $V$.
We will usually denote the Chow family associated to a family with the calligraphic version
of the same letter. 
\end{definition}

\begin{definition}
We denote by $\loc(\V^1, \dots, \V^k)$ the set of points $x \in X$ such that there exist
cycles $C_1, \dots, C_k$ with the following properties:
   \begin{itemize}
      \item $C_i$ belongs to the family $\V^i$;
      \item $C_i \cap C_{i+1} \neq \emptyset$;
      \item $x \in C_1 \cup \dots \cup C_k$,
   \end{itemize}
i.e. $\loc(\V^1, \dots, \V^k)$ is the set of points which belong to a connected chain of
$k$ cycles belonging \underline{respectively} to the families $\V^1, \dots, \V^k$.
\end{definition}

\begin{definition}
We denote by $\loc(\V^1, \dots, \V^k)_Y$ the set of points $x \in X$ such that there exist
cycles $C_1, \dots, C_k$ with the following properties:
   \begin{itemize}
      \item $C_i$ belongs to the family $\V^i$;
      \item $C_i \cap C_{i+1} \not = \emptyset$;
      \item $C_1 \cap Y \not = \emptyset$ and $x \in C_k$,
   \end{itemize}
i.e. $\loc(\V^1, \dots, \V^k)_Y$ is the set of points that can be joined to $Y$ by a connected
chain of $k$ cycles belonging \underline{respectively} to the families $\V^1, \dots, \V^k$.
\end{definition}

\begin{definition} 
Let $V^1, \dots, V^k$ be unsplit families on $X$.
We will say that $V^1, \dots, V^k$ are {\sf numerically independent} if their numerical classes
$[V^1], \dots ,[V^k]$ are linearly independent in the vector space $N_1(X)$.
If moreover $C \subset X$ is a curve we will say that $V^1, \dots, V^k$ are numerically 
independent from $C$ if the class of $C$ in $N_1(X)$ is not contained in the
vector subspace generated by $[V^1], \dots ,[V^k]$.   
\end{definition}

\begin{lemma} \label{locy}{\rm\cite[Lemma 5.4]{ACO}}
Let $Y \subset X$ be a closed subset and $V$ an unsplit family.
Assume that curves contained in $Y$ are numerically independent from curves in $V$, and that
$Y \cap \loc(V) \not= \emptyset$.
Then for a general $y \in Y \cap \loc(V)$
\begin{itemize}
      \item[\rm(a)] $\dim \loc(V)_Y \ge \dim (Y \cap \loc(V)) + \dim \loc(V_y);$
      \item[\rm(b)] $\dim \loc(V)_Y \ge \dim Y -K_X \cdot V - 1$.
\end{itemize}
Moreover, if $V^1, \dots, V^k$ are numerically independent
unsplit families such that curves contained in $G$ are numerically independent
from curves in $V^1, \dots, V^k$ then either 
$\loc(V^1, \ldots, V^k)_Y=\emptyset$ or
\begin{itemize}
      \item[\rm(c)] $\dim \loc(V^1, \ldots, V^k)_Y \ge \dim Y +\sum (-K_X \cdot V^i) -k$.
\end{itemize}
\end{lemma}

\begin{definition}
We define on $X$ a relation of {\sf rational connectedness with respect to $\V^1, \dots, \V^k$}
in the following way: $x$ and $y$ are in rc$(\V^1,\dots,\V^k)$-relation if there
exists a chain of rational curves in $\V^1, \dots ,\V^k$ which joins $x$ and $y,$ i.e.
if $y \in \cloc_m(\V^1, \dots, \V^k)_x$ for some $m$.
\end{definition}

To the rc$(\V^1,\dots,\V^k)$-relation we can associate a fibration, at least on an open subset.

\begin{theorem}{\rm\cite{Cam81}, \cite[IV.4.16]{Kob}} \label{rcvfibration}
There exist an open subvariety $X^0 \subset X$ and a proper morphism with connected fibers
$\pi\colon X^0 \to Z^0$ such that
   \begin{itemize}
      \item[\rm(a)] the rc$(\V^1,\dots,\V^k)$-relation restricts to an equivalence relation on $X^0;$
      \item[\rm(b)] the fibers of $\pi$ are equivalence classes for the rc$(\V^1,\dots,\V^k)$-relation$;$
      \item[\rm(c)] for every $z \in Z^0$ any two points in $\pi^{-1}(z)$ can be connected by a chain
      of at most $2^{\dim X - \dim Z}-1$ cycles in $\V^1, \dots, \V^k$.
   \end{itemize}
\end{theorem}

\begin{definition}
In the above assumptions, if $\pi$ is the constant map we say that $X$ is 
{\sf rc$(\V^1,\dots,\V^k)$-connected}.
\end{definition}

\begin{definition} 
A {\sf minimal horizontal dominating family}
with respect to $\pi$ is a family $V$ of horizontal curves such that $\loc(V)$
dominates $Z^0$ and $-K_X \cdot V$ is minimal among the families
with this property.\\
If $\pi$ is the identity map we say that $V$ is a {\sf minimal dominating family} for $X$.
\end{definition}

\begin{definition}
Let $\V$ be the Chow family associated to a family of rational curves $V$. We say that
$V$ is {\sf quasi-unsplit} if every component of any reducible cycle in $\V$ is
numerically proportional to $V$.\\
We say that $V$ is {\sf locally quasi-unsplit} if, for a general $x \in \loc(\V)$
every component of any reducible cycle in $\V_x$ is numerically proportional to $V$
\end{definition}

Note that any family of deformations of a rational curve whose numerical class lies in an extremal ray 
of $\cone(X)$ is quasi-unsplit.

\begin{lemma}\label{findunsplit} 
Let $X$ be a manifold and let $L$ be a line bundle on $X$ 
Let $V$ be a family of rational curves such that $L \cdot V > 0$.
Then there exists an unsplit family $V^L$ such that $L \cdot V^L >0$ and
$$[V] \equiv [V^L]+ [\Delta],$$
where $\Delta$ is an effective rational one cycle.
\end{lemma}

\begin{proof} 
If $V$ is unsplit there is nothing to prove, so assume that
the associated Chow family $\V$ contains a reducible cycle $\sum \Gamma_i$: then
for at least one $i$ we have $L \cdot \Gamma_i >0$. \\
Let $V^i$ be a family of deformations
of $\Gamma_i$; if $V^i$ is unsplit set $V^L=V^i$, otherwise let
$\sum \Gamma_{ij}$ be a reducible cycle in the associated Chow family $\V^i$:
then for at least one $j$ we have $L \cdot \Gamma_{ij} >0$.\\
Let $V^{ij}$ be a family of deformations
of $\Gamma_{ij}$; if $V^{ij}$ is unsplit set $V^L=V^{ij}$, otherwise 
continue as above. Since the degree of $V$ with respect to an ample line bundle is
finite the procedure ends after a finite number of steps.\end{proof}

\bigskip
{\bf Notation}: Let $S$ be a subset of $X$. 
We write $N_1(S)=\morespan{V^1}{V^k}$ if the numerical class in $\cycl(X)$ of every curve 
$C \subset S$ can be written as $[C]= \sum_i a_i [C_i]$, with $a_i \in \quadr$ and 
$C_i \in V^i$. We write $\cone(S)=\morespan{V^1}{V^k}$  (or $\cone(S)=\morespan{R_1}{R_k}$) 
if the numerical class in $\cycl(X)$ of every curve $C \subset S$ can
be written as $[C]= \sum_i a_i [C_i]$, with $a_i \in \quadr_{\ge 0}$ 
and $C_i \in V^i$ (or $[C_i]$ in $R_i$).\par

\ms
\begin{lemma} {\rm \cite[Lemma 1.4.5]{BSW}, \cite[Lemma 1]{Op}, \cite[Corollary 2.23]{CO}} \label{numequns}
Let $Y \subset X$ be a closed subset and $V$ an unsplit family of rational curves. 
Then every curve contained in $\loc(V)_Y$ is numerically equivalent to
a linear combination with rational coefficients
   $$a C_Y + b C_V,$$
where $C_Y$ is a curve in $Y,$ $C_V$ belongs to the family $V$ and $a \ge 0$.\\
Moreover, if $\Sigma$ is an extremal face of $\cone(X)$, $Y$ is a fiber
of the associated contraction and $[V]$ does not belong to $\Sigma,$
then
$$\cone(\cloc_m(V)_Y) = \langle \Sigma, [V] \rangle \quad \mathrm{for~every}\quad m \ge 1.$$
\end{lemma}


\section{Dominating families and Picard number}

We collect in this section some technical result that we will need in the proof.\\
The first is a variation of a classical construction of Mori theory, and says that, given a family
of rational curves $V$ and a curve $C$ contained in $\loc(V_x)$ for an $x$ such that $V_x$ is proper
we have $[C] \equiv a[V]$.\\
The only new remark - which already followed from the old proofs, but, to our best knowledge,
was not stated - is the fact that $a$ is a positive integer.

\begin{lemma}\label{intmult} 
Let $X$ be a smooth variety, $V$ a family of
rational curves on $X$, $x \in \loc(V)$ a point such that $V_x$ is proper
and $C$ a curve contained in $\loc(V_x)$.\\
Then $C$ is numerically equivalent to an integral multiple of a curve in $V$.
\end{lemma}

\begin{proof}  Consider the basic diagram
\begin{equation}
\xymatrix{p^{-1}(V_x)=:U_x  \ar[r]^(.65){i} \ar[d]_{p} & X\\
  V_x & & }
\end{equation}

Let $C$ be a curve contained in $\loc(V_x)$; if $C$ is a curve
parametrized by $V$ we have nothing to prove, so we can suppose that this is not the case.\\
In particular we have that $i^{-1}(C)$ contains  an irreducible curve $C'$ which is not
contained in a fiber of $p$ and dominates $C$ via $i$; let $S'$ be the surface $p^{-1}(p(C'))$,
let $B'$ be the curve $p(C') \subset V_x$ and let $\nu\colon B \to B'$ be the normalization of $B'$.
By base change we obtain the following diagram
$$
\xymatrix{S_B \ar[d]\ar[r]^{\bar \nu} & U_x \ar[d]^p \ar[r]^{i} &X\\
  B \ar[r]_\nu & V_x &}
$$
Let now $\mu\colon S \to S_B$ be the normalization of $S_B$; by standard arguments 
(see for instance \cite[1.14]{Wile})
it can be shown that $S$ is a ruled surface over the curve $B$; 
let $j\colon S \to X$ be the composition of $i$, $\bar \nu$ and $\mu$.
Since every curve parametrized by $S$ passes through $x$ there exists an irreducible 
curve $C_x \subset S$ which is contracted by $j$; by \cite[II.5.3.2]{Kob} we have
$C_x^2 <0$, hence $C_x$ is the minimal section of $S$.\\ 
Since every curve in $S$ is algebraically equivalent to a linear combination 
with integral coefficients of $C_x$ and a fiber $f$, and since
$C_x$ is contracted by $j$, every curve in $j(S)$ is algebraically equivalent in $X$
to an integral multiple of  $j_*(f)$, which is a curve of the family $V$;
but algebraic equivalence implies numerical equivalence and so the lemma is proved.\end{proof}

\begin{corollary}\label{kedu} 
Let $X$ be a smooth variety of dimension $n$ and let $V$ be a locally unsplit 
dominating family such that $-K_X \cdot V= n+1$; then $X \simeq \proj^n$.
\end{corollary}

\begin{proof} \quad For a general point  $x \in X$ we know that $V_x$ is proper and $X=\loc(V_x)$
by Proposition \ref{iowifam} (b). Therefore, by Lemma \ref{intmult}, for every curve $C$ in $X$
we have $-K_X \cdot C \ge n+1$ and we can apply \cite[Theorem 1.1]{Kepn}. \end{proof}

\begin{remark} The corollary also followed from the arguments in the proof
of \cite[Theorem 1.1]{Kepn}.
\end{remark}

In the rest of the section we establish some bounds on the Picard number of
Fano manifolds with minimal dominating families of high anticanonical degree.

\begin{lemma}\label{minimalD} 
Let $X$ be a Fano manifold of dimension $n \ge 3$ 
and pseudoindex $i_X \ge 2$ with a minimal dominating family $W$ such that $-K_{X} \cdot W > 2$;
if $X$ contains an effective divisor $D$ such that $\cone(D)=\onespan{[W]}$ then $\rho_{X}=1$.
\end{lemma}

\begin{proof} \quad The effective divisor $D$ has positive intersection number
with at least one of the extremal rays of $X$. Let $R$ be such a ray, denote
by $\f_R$ the associated contraction and by $V^R$ a family of deformations of a 
minimal rational curve in $R$.\\
If the numerical class of $W$ does not belong to $R$ then $D$ cannot contain
curves whose numerical class is in $R$, therefore every fiber
of $\f_R$ is one-dimensional.\\ 
By Proposition \ref{fiberlocus} this is possible only if 
$l(R) \le 2$ and therefore, since $l(R) \ge i_X$, it must be $l(R)=i_X = 2$.\\
Since every fiber of $\f_R$ is one-dimensional we have, for every $x \in \loc(V^R)$
that $\dim \loc(V^R_x)=1$ and therefore, by Proposition \ref{iowifam} (a) $V^R$ is a dominating
family. But, recalling that 
$$2 = -K_X \cdot V^R < -K_X \cdot W,$$
we contradict the assumption that $W$ is minimal.\\
It follows that $[W] \in R$, so the family $W$ is quasi-unsplit and $D \cdot W >0$;
hence $X$ can be written as $X=\loc(\cW)_D$, and 
by Lemma \ref{numequns} we have $\rho_X=1$. \end{proof}

\begin{corollary}\label{minimaln} 
Let $X$ be a Fano manifold of dimension $n \ge 3$ 
and pseudoindex $i_X \ge 2$ which admits a minimal dominating family $W$
such that $-K_{X} \cdot W \ge n$; then $\rho_{X}=1$.
\end{corollary}

\begin{proof} \quad 
Let $x \in X$ be a general point; every minimal dominating family
is locally unsplit, hence $\cone(\loc(W_x))=\onespan{[W]}$ by Lemma \ref{numequns}.\\
By Proposition \ref{iowifam} we have $\dim \loc(W_x) \ge -K_X \cdot W-1 \ge n-1$, so 
either $X= \loc(W_x)$ or $\loc(W_x)$ is an effective divisor
verifying the assumptions of Lemma \ref{minimalD}. 
In both cases we can conclude that $\rho_X=1$.   \end{proof}

\begin{lemma}\label{minimaln-1} 
Let $X$ be a Fano manifold of dimension $n \ge 3$ 
and pseudoindex $i_X \ge 2$, with a minimal dominating family $W$ such that
$-K_{X} \cdot W=n-1$; let $U \subset X$ be the open subset of 
points $x \in X$ such that $W_x$ is unsplit.
If a general curve $C$ of $W$ is contained in $U$ then either 
$\loc(W)_{C}$ is a divisor and $\rho_{X}=1$
or there exists an unsplit family $V$ such that $-K_X \cdot V=2$,
$D:=\loc(V)$ is a divisor and $D \cdot W >0$.
\end{lemma}

\begin{proof} \quad Let $C$ be a general curve in $W$ and consider
$\loc(W)_{C}$; by our assumptions we have $\cone(\loc(W)_{C})=\onespan{[W]}$
and $\dim \loc(W)_{C} \ge n -2$.\\
If $X=\loc(W)_{C}$ then clearly $\rho_X=1$, while if 
$\loc(W)_{C}$ has codimension one we conclude by Lemma \ref{minimalD}.\\
Therefore we can assume that, for a general $C$ in $W$, each component of 
$\loc(W)_{C}$ has codimension two in $X$. The fibration
$\pi\colon\rc{X}{Z}$ associated to the open prerelation defined by $W$ is proper, since a
general fiber $F$ coincides with $\loc(W_x)$ for a general $x \in F$ and 
$\loc(W_x)$ is closed since $W$ is locally unsplit.\\
Being $\pi$ proper there exists a  minimal horizontal dominating family $V$
with respect to $\pi$; since the general fiber of $\pi$ has dimension
$n-2$, then $\dim Z=2$, hence for a general $x \in \loc(V)$ we have
$\dim \loc(V_x) \le 2$.\\
It follows that $V$ is an unsplit family, which cannot be dominating by the minimality of $W$, so 
$\dim \loc(V_x) \ge i_{X} \ge 2$, and $D=\loc(V)$ is a divisor by Proposition \ref{iowifam}.
Since $D$ dominates $Z$ we have $D \cdot W > 0$.\end{proof}


\section{Fano manifolds obtained blowing-up non Fano manifolds}

We start now the proof of our results. Let us fix once and for all the 
setup and the notation:

\begin{statement} \label{setup}
{\it $X$ is a Fano manifold of pseudoindex $i_X \ge 2$ and dimension $n \ge 6$, 
which has a contraction $\sigma\colon X \to Y$ which is the blow-up of a manifold 
$Y$ along a smooth subvariety $B$ of dimension $i_X$. We denote by $R_\sigma$
the extremal ray corresponding to $\sigma$, by $l_\sigma$ its length and
by $E$ its exceptional locus.} 
\end{statement}

\begin{remark} The assumption on $\dim B$ is equivalent to
$$l_\sigma +i_X = n-1.$$
\end{remark}

In this section we will deal with Fano manifolds as in Theorem \ref{main}
which are obtained as a blow-up $\sigma \colon X \to Y$ of a manifold $Y$
which is not Fano. It turns out that there is only one possibility 
(Corollary \ref{crelle}) we start with a slightly general result:

\begin{theorem}\label{Eneg1} 
Let $X$, $R_\sigma$ and $E$ be as in \ref{setup} and assume that there exists on 
$X$ an unsplit family of rational curves $V$ such that $E \cdot V  <0$; then either 
$[V] \in R_\sigma$ or $X =\proj_{\proj^{n-3} \times \pd}(\Ol(1,1) \oplus \Ol(2,2))$.
\end{theorem}

\begin{proof} 
Since $E \cdot V <0$ then $\loc(V) \subseteq E$, so $V$ is not a dominating
family.\\
Pick $x \in \loc(V)$ and let $F_\sigma$ be the fiber of $\sigma$ through $x$; we have 
$$\dim E \ge \dim \loc(V_x)+ \dim F_\sigma \ge i_X +l_\sigma \ge n-1,$$
so all the above inequalities are equalities; in particular we have $\dim \loc(V_x)=i_X$
and so, by Proposition \ref{iowifam},
$$\dim \loc(V) \ge n + i_X-1- \dim \loc(V_x) = n-1,$$
hence $\loc(V)=E$; therefore the above (in)equalities are true for every $x \in E$.
It follows that $\sigma$ is equidimensional and so it is a smooth blow-up by 
\cite[Theorem 5.1]{AOsperays}.\par

\smallskip
Considering $V$ as a family on the smooth variety $E$ we can write
$$n-1+i_X=\dim \loc(V) + \dim \loc(V_x) \ge -K_E \cdot V + n-2,$$
therefore $-K_E \cdot V \le i_X +1$; on the other hand 
$$-K_E \cdot V = -K_X \cdot V -E \cdot V \ge i_X +1,$$
forcing $-K_E \cdot V= i_X +1$ and $E \cdot V=-1$.\\
Then on $E$ we have two unsplit dominating families of rational curves verifying
the assumptions of \cite[Theorem 1]{Op}, hence $E \simeq \proj^{i_X} \times \proj^{l_\sigma}$; 
in particular $\rho_E=2$.\\
Now let $R$ be an extremal ray of $X$ such that $E \cdot R >0$; by \cite[Corollary 2.15]{NO}
the contraction $\f_R$ associated to $R$ is a $\pu$-bundle; in particular,
by Proposition \ref{fiberlocus}, this implies that $i_X=2$.\\ 
Moreover, denoted by $V^R$ a family of deformation of a minimal rational curve
in $R$, we have $X=\loc(V^{R})_E$, so $\rho_X=3$ and
the description of $X$ is obtained arguing as in the 
proof of Proposition 7.3 in \cite{NO}.\end{proof}

\begin{corollary} \label{crelle}
In the assumptions of Theorem \ref{main} either $Y$ is a Fano manifold or 
$X =\proj_{\proj^{n-3} \times \pd}(\Ol(1,1) \oplus \Ol(2,2))$, 
$Y \simeq \proj_\pd(\Ol \oplus \Ol(1)^{n-2})$
and $B \simeq \pd$ is the section corresponding to the surjection $\Ol \oplus \Ol(1)^{n-2} \to \Ol$.
\end{corollary}

\begin{proof} If $Y$ is not Fano then by \cite[Proposition 3.4]{Wicon} there exists 
an extremal ray $R' \in \cone(X)$ such that $E \cdot R' <0$. \end{proof}

\begin{remark} Note that, if $X \simeq \proj_{\proj^{n-3} \times \pd}(\Ol(1,1) \oplus \Ol(2,2))$,
then $\cone(X)$ is generated by three extremal rays: one -- the $\pu$-bundle contraction -- is 
of fiber type, while the other two are smooth blow-ups with the same exceptional locus.
In particular $\cone(X)$ is as in Theorem \ref{main}, case (d).
\end{remark}

\begin{corollary}\label{Eneg2}
Let $X$, $R_\sigma$ and $E$ be as in \ref{setup}; assume that $Y$ is a Fano manifold
and that there exists on $X$ a family of rational curves $V$ such that $E \cdot V  <0$; then 
$-K_X \cdot V \ge l_\sigma$; moreover, if $V$ is unsplit then $[V] \in R_\sigma$.
\end{corollary}

\begin{proof} Arguing as in Lemma \ref{findunsplit} we can find an unsplit family $V^E$ such that 
$E \cdot V^E <0$ and $[V] \equiv [V^E]+ [\Delta]$. By Theorem \ref{Eneg1} we have that
$[V^E] \in R_\sigma$, hence 
$$-K_X \cdot V \ge -K_X \cdot V^E \ge l_\sigma.$$
To prove the last assertion note that, if $V$ is an unsplit family, we can apply 
Theorem \ref{Eneg1} directly to $V$. \end{proof}


\section{Manifolds with a dominating (quasi)-unsplit family} \label{unsplitsec}

In this section we will describe the cone of curves of Fano manifolds
as in \ref{setup} which admit a dominating quasi-unsplit family of rational curves $W$,
and such that the target of the blow-up $\sigma\colon X \to Y$ is a Fano manifold.\par
\medskip
If the family $W$ is quasi-unsplit but not unsplit then the result can be obtained easily:

\begin{lemma}\label{quasi}
Assume that $W$ is not unsplit; then $\rho_X = 2$, $i_X=2$ and $\cone(X)=\twospan{R_\sigma}{[W]}$.
\end{lemma}

\begin{proof} Since $W$ is not unsplit we have $-K_X \cdot W \ge 2i_X$.
Consider the associated Chow family $\cW$ and the  rc$\cW$-fibration
$\pi\colon \rc{X}{Z}$; since a general fiber of $\pi$ contains $\loc(W_x)$ for some $x$,
and by Proposition \ref{iowifam} we have $\dim \loc(W_x) \ge -K_X \cdot W \ge 2i_X-1$ we have
$$\dim Z \le n+1 -2i_X \le n-1 -i_X= \dim F_\sigma,$$
where $F_\sigma$ is a fiber of $\sigma$.\\
A family $V^\sigma$ of deformations of a minimal curve in $R_\sigma$
is thereby horizontal and dominating with respect to $\pi$; moreover,
since $F_\sigma$ dominates $Z$ we have that $X=\loc(\cW)_{F_\sigma}$,
hence $\cone(X)=\twospan{R_\sigma}{[W]}$ by Lemma \ref{numequns}.\end{proof}

\medskip
In view of Lemma \ref{quasi}, we can assume throughout the section that
$W$ is an unsplit dominating family. 

\begin{lemma}\label{mortadella} 
Let $X$ be a Fano manifold with $\rho_X=3$. 
Assume that there exists an effective divisor
$E$ which is negative on one extremal ray $R$ of $\cone(X)$ and is nonnegative on the other extremal rays. 
If $E \cdot C=0$ for a curve $C \subset  X$ whose numerical class lies in $\partial \cone(X)$, 
then $[C]$ is contained in a two-dimensional face of $\cone(X)$ which contains $R$.
\end{lemma}

\begin{proof} \quad By assumption, neither $E$ nor $-E$ are nef, hence the hyperplane $\{E=0\}$ 
has nonempty intersection with the interior of $\cone(X)$.
Let $\Sigma$ be a two-dimensional 
face of $\cone(X)$ containing $[C]$: by the above discussion
$E$ cannot be trivial on the whole face $\Sigma$.\\
Therefore, if $[C]$ lies in the interior of $\Sigma$ then $E$ must be negative on one of the rays
spanning $\Sigma$, hence $R \in \Sigma$.
If $[C]$ lies on an extremal ray, then $E$ has different sign on the rays 
which span with $[C]$ a two-dimensional face of $\cone(X)$, so $E$ is negative on one of them, 
which has to be $R$.\end{proof}

\begin{lemma}\label{nosmall}
Assume that there exists an extremal ray $R_\tau$ such that $[W] \not \in R_\tau$ and
either $E \cdot R_\tau > 0$ or $E \cdot W >0$. 
Then every fiber of the contraction $\tau$ associated to $R_\tau$ 
has dimension not greater than two. 
In particular $\tau$ is either a fiber type contraction or a smooth blow-up of a codimension three subvariety,
and in this case the exceptional locus of $\tau$ is $\Exc(\tau)=\loc(W,V^\tau)_{F_\sigma}$, for some
fiber $F_\sigma$ of $\sigma$.
\end{lemma}

\begin{proof} \quad Let $F_\tau$ be a fiber of $\tau$. 
If $E \cdot R_\tau >0$ there exists a  fiber $F_\sigma$ of $\sigma$ which meets $F_\tau$; 
since $W$ is dominating we have $F_\sigma \subset \loc(W)_{F_\sigma}$ 
and therefore $F_\tau \cap \loc(W)_{F_\sigma} \not = \emptyset$.\\
If else $E \cdot W > 0$ then $E \cap \loc(W)_{F_\tau} \neq \emptyset$, so there exists a fiber
$F_\sigma$ of $\sigma$ such that $F_\sigma \cap \loc(W)_{F_\tau} \neq \emptyset$;
equivalently, we have that $F_\tau \cap \loc(W)_{F_\sigma} \neq \emptyset$.\\
In both cases, this intersection cannot be of positive dimension, since every curve
in $F_\tau$ has numerical class belonging to $R_\tau$, while every curve in $\loc(W)_{F_\sigma}$
has numerical class contained in the cone $\langle R_\sigma,[W] \rangle$.
By our assumptions
$$\dim \loc(W)_{F_\sigma} \ge \dim F_\sigma + i_X -1 \ge l_\sigma + i_X -1 \ge n-2,$$ 
hence $\dim F_\tau \le 2$. Proposition \ref{fiberlocus} implies that
$\tau$ cannot be a small contraction; if it is divisorial,
by the same inequality it is equidimensional with two-dimensional fibers,
so it is a smooth blow-up by \cite[Theorem 5.1]{AOsperays}.\\
In this last case, denoted by $V^\tau$ a family of deformations
of a minimal curve in $R_\tau$, we have 
$$\dim \loc(W,V^\tau)_{F_\sigma} \ge n-1,$$
hence $\Exc(\tau)=\loc(W,V^\tau)_{F_\sigma}$. \end{proof}

\begin{lemma}\label{timoteo}
Assume that $E \cdot W=0$. Let $\pi\colon \rc{X}{Z}$ be the rc$W$-fibration
and let $V$ be a minimal horizontal dominating family with respect to $\pi$.
Then $R_\sigma$, $W$ and $V$ are numerically independent. In particular
$\rho_X \ge 3$. 
\end{lemma}

\begin{proof}  Since $E \cdot W=0$, $E$ does not dominate $Z$, hence 
$E$ cannot contain $\loc(V)$ and therefore $E \cdot V  \ge 0$.\\
Let $\ac$ be the pull-back to $X$ of a very ample divisor in $\pic(Z)$;
$\ac$ is zero on curves in the family $W$ and it is positive outside the indeterminacy 
locus of $\pi$; in particular $\ac \cdot V > 0$.\\
If $[V]$ were contained in the plane spanned by $R_\sigma$ and $[W]$
we could write $[V]=\alpha [V^\sigma]+ \beta [W]$, but
intersecting with $E$ we would get $\alpha \le 0$, while intersecting with $\ac$
we would get $\alpha >0$, a contradiction which proves the lemma.\end{proof}

\begin{proposition}\label{tito}
Assume  that $E \cdot W=0$. 
Let $\pi$ be the rc$W$-fibration and let $V$ be a minimal horizontal 
dominating family with respect to $\pi$. Then $V$ is unsplit.
\end{proposition}

\begin{proof}  Assume first that $E \cdot V >0$.\\ 
If $V$ is not unsplit we will have, for a general $x \in \loc(V)$, that 
$$\dim \loc(V_x)  \ge 2 i_X -1 \ge 3.$$
Since $E \cdot V >0$, then $E \cap \loc(V_x) \not = \emptyset$, therefore
$\loc(V_x)$ meets a fiber $F_\sigma$  of $\sigma$. 
Moreover, since $W$ is dominating, $F_\sigma \subset \loc(W)_{F_\sigma}$ and so
the intersection $\loc(V_x) \cap  \loc(W)_{F_\sigma}$ is not empty. 
This fact, together with
$$\dim \loc(W)_{F_\sigma} \ge l_\sigma +i_X -1 \ge n-2,$$
implies that $\loc(W)_{F_\sigma}$ contains a curve whose class is proportional to $[V]$,
a contradiction by Lemma \ref{timoteo},
since $\cone(\loc(W)_{F_\sigma})=\twospan{[W]}{R_\sigma}$.\par
\medskip
We will now deal with the harder case $E \cdot V=0$, assuming by contradiction that 
$V$ is not unsplit.\par
\medskip
We claim that $E$ has non zero intersection number with at least one component 
of a cycle in the Chow family $\V$. 
To prove the claim, consider the rc$(W,\V)$-fibration $\pi_{W,\V}$;
a general fiber of $\pi_{W,\V}$  contains $\loc(V,W)_x$ for some $x$,
so it has dimension $\ge 3i_X -2$.\\
Since $E$ is not contained in the indeterminacy locus of $\pi_{W,\V}$ - which has codimension at
least two in $X$ - it meets some fiber $G$ of $\pi_{W,\V}$ which, by semicontinuity, has dimension 
$\ge 3i_X -2$. 
Therefore there exists a fiber $F_\sigma$ of $\sigma$ such that $F_\sigma \cap G \not = \emptyset$.
and, for such a fiber we have
$$\dim (F_\sigma \cap G) \ge l_\sigma +3i_X-2 -n \ge 2i_X -3 \ge 1;$$
Let $C$ be a curve in $F_\sigma \cap G$; since $C \subset F_\sigma$ we have
$E \cdot C <0$; on the other hand, since $C \subset G$ the numerical class of $C$ can be written
as a linear combination of $[W]$ and of classes of irreducible components of cycles in $\V$.
Since $E \cdot W=0$ we see that $E$ cannot have zero intersection number with all
the components of cycles in $\V$ and the claim is proved.\par
\medskip
So in $\V$ there exists a reducible cycle $\Gamma=\sum_{i=1}^k \Gamma_i$ such that $E \cdot \Gamma_1 <0$.
Applying Lemma \ref{findunsplit} we find an unsplit family $T$ on which $E$ is negative
and such that $[\Gamma_1]=[T]+[\Delta]$, with $\Delta$ an effective rational 1-cycle.
\par
Since $Y$ is a Fano manifold, by Corollary \ref{Eneg2} we have that 
$[T] \in R_\sigma$ and $-K_X \cdot T \ge l_\sigma$; therefore, for a general $x \in \loc(V)$
$$\dim \loc(V_x) \ge -K_X \cdot V-1 =-K_X \cdot (T + \Delta + \sum_{i=2}^k \Gamma_i)-1
\ge l_\sigma+i_X-1 \ge n-2.$$
If $\dim \loc(V_x) \ge n-1$ then $X=\loc(W)_{\loc(V_x)}$ and $\rho_X=2$ against Lemma \ref{timoteo};
therefore $\dim \loc(V_x)= -K_X \cdot V -1 = n-2$, hence $V$ is a dominating
family by Proposition \ref{iowifam}, $\Gamma = \Gamma_1 +\Gamma_2$, $\Delta =0$, 
$\Gamma_1 \in R_\sigma$ and $-K_X \cdot \Gamma_2 =i_X$.\par
\medskip
For a general $x \in \loc(V)$ we have $\dim \loc(W)_{\loc(V_x)} \ge n-1$ by Lemma \ref{locy};
moreover, since $\cone(D)=\twospan{[W]}{[V]}$ and $\rho_X \ge 3$ by Lemma \ref{timoteo}, we cannot
have $D = X$, hence $D$ is an effective divisor.\\
We will now reach a contradiction by showing that $D$ has zero intersection number
with every extremal ray of $X$.\\
Let $\overline{V}$ be any unsplit 
family whose numerical class is not contained in
the plane spanned by $[W]$ and $[V]$; we cannot have
$\dim \loc(\overline{V}_x)=1$, otherwise $\overline{V}$ would be dominating of 
anticanonical degree $2$, against the minimality of $V$.
This implies that $D \cdot \overline{V}=0$ since $\cone(D)=\twospan{[W]}{[V]}$ implies
that $D \cap \loc(\overline{V}_x) = \emptyset$.\\
It follows that $D \cdot \Gamma_2 =0$ and that $D$ is tri\-vial on every extremal ray not lying
in the plane $\twospan{[V]}{[W]}$.
Since $[V]=[\Gamma_1]+[\Gamma_2]$ and $\Gamma_1 \in R_\sigma$, which is a ray
not contained in the plane spanned by $[W]$ and $[V]$ we have
that also  $D \cdot V =0$.\\
To conclude it is now  enough to observe that we must have $D \cdot W=0$, otherwise $\cloc_2(W)_{\loc(V_x)} = X$,
forcing again $\rho_X=2$. We have thus reached a contradiciton, since the effective
divisor $D$ has to be trivial on the whole $\cone (X)$.\end{proof}

\begin{proposition}\label{replace} Up to replace $W$ with another 
dominating unsplit family, we can assume that $E \cdot W > 0$.
\end{proposition}

\begin{proof}Assume that $E \cdot W=0$, let $\pi$ be the rc$W$-fibration, and let $V$ be 
a minimal horizontal dominating 
family with respect to $\pi$. By Proposition \ref{tito} we know that $V$ is unsplit.\par
\medskip
{\bf Case a)} \quad $V$ is dominating.\par
\medskip
If $E \cdot V > 0$ the Proposition is proved, so we can assume that $E \cdot V=0$.\\
If $F_\sigma$ is any fiber of $\sigma$ we have 
$$\dim \loc(V,W)_{F_\sigma} \ge \dim F_\sigma +2i_X-2 = l_\sigma +2 i_X -2\ge n-1.$$
Note that, by the assumptions on the intersection numbers, we have $\loc(V,W)_{F_\sigma} \subseteq E$,
and therefore $\loc(V,W)_{F_\sigma} = E$; in particular it follows from 
the above inequalities that $i_X=2$.\\
We can repeat the same arguments to show that also $\loc(W,V)_{F_\sigma} = E$; 
hence every curve contained in $E$ is numerically equivalent to a linear combination
$$a [V^\sigma] + b [V] + c [W]$$ 
with $a$, $b$, $c \ge 0$ by Lemma \ref{numequns},
and therefore
$\cone(E)=\threespan{R_\sigma}{[V]}{[W]}$. In particular 
$E$ has nonpositive intersection with every curve it contains.\\
Let $R_\vartheta$ be an extremal ray such that $E \cdot R_\vartheta>0$; by
\cite[Corollary 2.15]{NO} the associated contraction $\vartheta\colon X \to Y$
is a $\pu$-bundle; the associated family $V^\vartheta$ is dominating and unsplit
and $E \cdot V^\vartheta >0$, and the proposition is proved.\par
\medskip
{\bf Case b)} \quad $V$ is not dominating.\par
\medskip
Consider the rc$(W, V)$-fibration $\pi'\colon \rc{X}{Z'}$; $Z'$ has positive dimension since
by Lemma \ref{timoteo} we have $\rho_X \ge 3$.\\
A general fiber $F'$ of $\pi'$ contains $\loc(V,W)_x$ for some $x \in \loc(V)$, hence
$\dim F' \ge 2i_X -1$ and thus 
$$\dim Z' \le n+1 -2i_X \le l_\sigma.$$
A general fiber $F_\sigma$ of $\sigma$ is not contained in the indeterminacy locus of $\pi'$ and 
is not contracted by $\pi'$,
since, by Lemma \ref{timoteo}, $[V]$, $[W]$ and  $R_\sigma$ are numerically independent. 
Hence we have $\dim Z' \ge \dim F_\sigma = l_\sigma$ and the above inequalities
are equalities.\\ 
It follows that $i_X=2$, $\dim Z' = l_\sigma$ and $F_\sigma$
dominates $Z'$; this implies that $X = \cloc_m(W,V)_{F_\sigma}$ for some $m$.
Therefore $\rho_X=3$ and so, by Lemma \ref{numequns}, the numerical class of every curve in 
$X$ can be written as 
$$\alpha[V^\sigma]+\beta[W] +\gamma[V],$$
with $\alpha \ge 0$. This implies that the plane $\langle [V],[W] \rangle$ is extremal in $\cone(X)$.\par

\medskip
The divisor $E$ has to be positive on $V$, otherwise it would be nonpositive on the whole $\cone(X)$;
since $E \cdot W=0$ then $[W]$ is in an extremal face with $R_\sigma$ by Lemma \ref{mortadella}.
Since $[W]$ is also in an extremal face with $[V]$ it follows that $[W]$
spans an extremal ray of $\cone(X)$, whose associated contraction is of fiber type.\\

Let $W_Y$ be a minimal dominating family on $Y$ and let $W^*$ be a family of deformations of the 
strict transform of a general curve in $W_Y$.\\
We have $-K_X \cdot W^* =-K_Y \cdot W_Y \le n$; in fact, if   
$-K_X \cdot W^* = n+1$, we would have $Y \simeq \proj^n$ by Corollary \ref{kedu} 
and so $\rho_X=2$.\\
Assume that $W^*$ is not locally unsplit; then there
exists a reducible cycle $\Gamma=\sum \Gamma_i$ in $\cW^*$ such that the family $T$
of deformation of one irreducible component, say $\Gamma_1$ is dominating.\\
We cannot have $E \cdot T = 0$, otherwise, denoting by $T_*$ a family of deformation
of the image in $Y$ of a general curve in $T$ we would have
$$-K_X \cdot T = -K_{Y} \cdot T_* < -K_{Y} \cdot W_Y$$
against the minimality of $W_Y$.\\
Therefore $E \cdot T >0$; in this case $E$ must be negative on another
component of $\Gamma$, say $\Gamma_2$.
By Corollary \ref{Eneg2}, we have that $-K_X \cdot \Gamma_2 \ge l_\sigma$, and thus
$\Gamma = \Gamma_1 + \Gamma_2$ and so $-K_X \cdot \Gamma_1 < 2i_X$ and the family $T$
is unsplit and dominating.\\
We are left with the case in which $W^*$ is locally unsplit.\\
Consider $\loc(V^\sigma)_{\loc(W^*_x)}$ for a general $x \in \loc(W^*)$: it is contained in $E$
and, by Lemma \ref{locy} 
$$\dim(\loc(V^\sigma)_{\loc(W^*_x)}) \ge n-2+l_\sigma-1,$$
yielding $l_\sigma \le 2$ and so $n=5$, a contradiction which concludes the proof.
\end{proof}

\begin{theorem}\label{part1}
Let $X$ be a Fano manifold of pseudoindex $i_X \ge 2$ and dimension $n \ge 6$, with a contraction
$\sigma\colon X \to Y$ which is the blow-up of a Fano manifold $Y$ along a smooth subvariety 
$B$ of dimension $i_X$.
If $X$ admits a dominating unsplit family of rational curves $W$
then the possible cones of curves of $X$ are listed in the following table,
where $R_\sigma$ is the ray corresponding to $\sigma$, $F$ stands for a fiber type 
contraction and $D_{n-3}$ for a divisorial contraction whose
exceptional locus is mapped to a subvariety of codimension three.
\begin{center}
\begin{tabular}{|c|c||c|c|c|c|c|}\hline
\quad $\rho_X$ \quad & \quad $i_X$ \quad &  \quad $R_1$
\quad & \quad $R_2$ \quad & \quad $R_3$ \quad & \quad $R_4$ \quad \\
\hline\hline
  2 &   & $R_\sigma$  & $F$ &  &  \\
  3 & 2,3 & $R_\sigma$  & $F$ & $F$ &   \\
  3 & 2 & $R_\sigma$ & $F$ & $D_{n-3}$ &  \\
  4 &  2   & $R_\sigma$ & $F$  & $F$ & $F$\\
\hline
\end{tabular}
\end{center}~\par
In particular generalized Mukai conjecture holds for $X$.
\end{theorem}

\begin{proof}
Let $V^\sigma$ be a family of deformations of a minimal rational curve in $R_\sigma$.\\
By Proposition \ref{replace} we can assume that $E \cdot W >0$; therefore
the family $V^\sigma$ is horizontal and dominating with respect
to the rc$W$-fibration $\pi\colon\rc{X}{Z}$.\\
It follows that a general fiber $F'$ of the the rc$(W,V^\sigma)$-fibration $\pi'\colon\rc{X}{Z'}$
contains $\loc(W)_{F_\sigma}$ for some fiber $F_\sigma$ of $\sigma$, and therefore
$$ \dim F' \ge \dim \loc(W)_{F_\sigma} \ge l_\sigma+i_X-1 \ge n-2,$$ 
hence $\dim Z' \le 2$.\par
\smallskip
If $\dim Z'=0$ then $X$ is rc$(W,V^\sigma)$-connected and $\rho_X=2$; denote by
$R_\vartheta$ the extremal ray of $\cone(X)$ different from $R_\sigma$.
We claim that in this case $[W] \in R_\vartheta$. In fact, if this were not the case, 
$R_\vartheta$ would be a small ray by \cite[Lemma 2.4]{CO}, but in our assumptions 
we have $E \cdot R_\vartheta >0$, against Lemma \ref{nosmall}.\\
We can thus conclude that in this case
$\cone(X) = \twospan{R_\sigma}{R_\vartheta}$ and that $R_\vartheta$ is of fiber type.\par
\medskip
If $\dim Z' > 0$ take $V'$ to be a minimal horizontal dominating family for $\pi'$; by
\cite[Lemma 6.5]{ACO} we have $\dim \loc(V'_x) \le 2$, and therefore
$$-K_X \cdot V' \le \dim \loc(V'_x) +1 \le 3,$$
so $V'$ is unsplit and $i_X \le 3$.\par 
\medskip
The classes $[V^\sigma]$ and $[W]$ lie on an extremal face $\Sigma=\twospan{R_\sigma}{R}$ 
of $\cone(X)$, since, otherwise, by \cite[Lemma 2.4]{CO}, 
$X$ would have a small contraction, against Lemma \ref{nosmall}.
Let $\ac$ the pull back via $\pi$ of a very ample divisor on $Z$.\\
We know that $\ac \cdot W = 0$ and $\ac \cdot R_\sigma >0$, since 
$V^\sigma$ is horizontal and dominating with respect to $\pi$.
It follows that $[W] \in R$ (and so $R$ is of fiber type), since otherwise
the exceptional locus of $R$ would be contained in the indeterminacy locus of $\pi$,
and thus the associated contraction would be small, contradicting again
Lemma \ref{nosmall}.\par
\medskip
Consider now the rc$(W,V^\sigma,V')$-fibration $\pi''\colon\rc{X}{Z''}$: 
its fibers have dimension $\ge n-1$ and so $\dim Z'' \le 1$.\par
\medskip
If $\dim Z''=0$ we have that $X$ is rc$(W,V^\sigma,V')$-connected and $\rho_X=3$; 
by Lemma \ref{nosmall} every extremal ray of $X$ has an associated contrsction
which is either of fiber type or divisorial.\par
\smallskip
Assume that there exists an extremal ray $R'$ not belonging to $\sigma$ such that
its associated contraction is of fiber type.\\
This ray must lie in a face of $\cone(X)$ with $R$ by \cite[Lemma 5.4]{CO}.\\ 
If $E \cdot R' > 0$ we can exchange the role of $R$ and $R'$ and repeat the previous argument,
therefore $R'$ lies in a face with $R_\sigma$ and $\cone(X)=\threespan{R_\sigma}{R}{R'}$. \\
If $E \cdot R' = 0$ there cannot be any extremal ray 
in the half-space of $\cone(X)$ determined by the plane 
$\twospan{R'}{R_\sigma}$ and not containing $R$, otherwise this ray would have negative 
intersection with $E$, a contradiction. So again $\cone(X)=\threespan{R_\sigma}{R}{R'}$.\par
\medskip
So we can assume that every ray not belonging to $\Sigma$ is divisorial.
Let $R'$ be such a ray, denote by $E'$ its exceptional locus, and by $W'$ a family 
of deformations of a minimal rational curve in $R'$.\\
Let $F'$ be a fiber of the rc$(W,V^\sigma)$-fibration $\pi'$;
since $\dim F' \ge n-2$ we can write $E'=\loc(W')_{F'}$.
It follows that $\cone(E')=\threespan{R_\sigma}{R}{R'}$.
In particular $E'$ cannot be trivial on $\Sigma$, otherwise it would be nonpositive 
on the whole $\cone(X)$.\\
We claim that $R$ and $R'$ lie on an extremal face of $\cone(X)$:
if $E' \cdot R > 0$ the family $W'$ is horizontal and dominating with respect to $\pi$ and so
$R'$ and $R$ are in a face by \cite[Lemma 5.4]{CO}.
If else $E' \cdot R = 0$ we have $E' \cdot R_\sigma > 0$. It follows that 
in the half-space determined by $\twospan{R}{R'}$ 
and not containing $R_\sigma$ the divisor $E'$ is negative.
Therefore this half space cannot contain an extremal ray $R''$, since otherwise, 
the exceptional locus of this ray must be contained in $E'$, contradicting the fact 
that $\cone(E')=\threespan{R_\sigma}{R}{R'}$.\\
So we have proved that every ray not belonging to $\Sigma$ lies in a face with $R$, 
and this implies that such a ray is unique and $\cone(X)=\threespan{R_\sigma}{R}{R'}$.\par
\medskip
Recalling that $E'=\loc(W')_{F'}$ and that $\dim F' \ge n-2$ we have that every fiber of 
the contraction $\f'$ associated to $R'$ has dimension two; it follows that $i_X=2$ and that
$\f'$ is a smooth blow-up of a codimension three subvariety by \cite[Theorem 5.1]{AOsperays}.\par
\bigskip
If $\dim Z'' =1$  consider a minimal horizontal dominating family $V''$ for $\pi''$: in this case 
$\rho_X=4$, $i_X=2$ and  both $V'$ and $V''$ are dominating.
Let $F_\sigma$ be a fiber of $\sigma$: then we can write 
$X = \loc(V', V'')_{\loc(W)_{F_\sigma}}$. By Lemma \ref{numequns} every curve in $X$ can 
be written with positive coefficients with respect to $V^\sigma$ and $W$; but $W$, $V'$ and $V''$ 
play a  symmetric role, so we can conclude
that $\cone(X) = \langle R_\sigma, [W], [V'], [V''] \rangle$, and all the three 
rays different from $R_\sigma$
are of fiber type. \end{proof}

\section{Manifolds without a dominating quasi-unsplit family}

In this section we will show that the only Fano manifold
as in \ref{setup} which does not admit a dominating quasi-unsplit family of rational curves
is the blow-up of $\G(1,4)$ along a plane of bidegree $(0,1)$ (Theorem \ref{noqu}).
In view of Theorem \ref{part1} this will conclude the proof
of Theorem \ref{main} and prove Theorem \ref{grass}.\par 
\medskip

From now on we will thus work in the following setup:

\begin{statement} \label{setup2}
{\it $X$ is a Fano manifold of pseudoindex $i_X \ge 2$ and dimension $n \ge 6$, 
which does not admit a quasi-unsplit dominating family of rational curves and 
has a contraction $\sigma\colon X \to Y$ which is the blow-up of a manifold 
$Y$ along a smooth subvariety $B$ of dimension $i_X$. We denote by $R_\sigma$
the extremal ray corresponding to $\sigma$, by $l_\sigma$ its length and
by $E$ its exceptional locus.} 
\end{statement}

In view of Corollary \ref{crelle} we can assume that $Y$ is a Fano manifold.
We need some preliminary work to establish some properties of families
of rational curves on $X$ and $Y$.

\begin{lemma}\label{doccia} 
Assume that $\rho_X=2$. Let $W'$ be a minimal dominating family of rational curves for $Y$
and let $W^*$ be a family of deformations of the strict transform of a general curve in $W'$.
Then $-K_{Y} \cdot W' \ge n-1$.
\end{lemma}

\begin{proof} \quad The family $W^*$ is dominating and therefore, by \ref{setup2}, 
not quasi unsplit. Moreover, by Corollary \ref{E0} we have $E \cdot W^*=0$, hence
there exists a component $\Gamma^*_1$ of a reducible cycle 
$\Gamma^*$ in $\cW^*$ such that $E \cdot \Gamma^*_1 <0$.\\
By Corollary \ref{Eneg2} we have $-K_X \cdot \Gamma^*_1 \ge l_\sigma$, and therefore
 $$-K_{Y} \cdot W' = -K_X \cdot W^* \ge l_\sigma+i_X = n-1. $$
  \end{proof}

\begin{proposition} \label{p5}
Let $X$, $Y$, $R_\sigma$ and $E$ be as in \ref{setup2}. 
Then there does not exist on $X$ any locally unsplit dominating family $W$ 
such that $E \cdot W >0$. 
\end{proposition}

\begin{proof} \quad Assume  that such a family $W$ exists; we will derive a contradiction showing that
in this case $n = 5$.\par

\medskip
First of all we prove that $i_X=2$ and that $X$ is rationally connected
with respect to the Chow family $\cW$ associated to $W$ and to $V^\sigma$, the family of 
deformations of a general curve of minimal degree in $R_\sigma$.\\
Since $E \cdot W >0$, for a general $x \in X$, the intersection 
$E \cap \loc(W_x)$ is nonempty. 
On the other hand, the fact that $E \cdot V^\sigma <0$ yields  
that the families $W$ and $V^\sigma$ are numerically independent, and therefore,
for every fiber $F_\sigma$ of $\sigma$, we have $\dim(\loc(W_x) \cap F_\sigma) \le 0$.\\ 
Now, if we denote by $F_\sigma$ a fiber of $\sigma$
which meets $\loc(W_x)$, it follows that
$$2i_X -1 \le -K_X \cdot W -1 \le \dim \loc(W_x) \le n - \dim F_\sigma \le n - l_\sigma =i_X+1,$$ 
whence $i_X=2$, $\dim \loc(W_x) = i_X+1 =3$ and $-K_X \cdot W = 2i_X = 4$.\\
In particular
$\dim (E \cap \loc(W_x)) = 2 = \dim B$, hence $\sigma(E \cap \loc(W_x))=B$
and every fiber of $\sigma$ meets $\loc(W_x)$.\\
Let $x$ and $y$ be two general points in $X$; every fiber of $\sigma$ meets both $\loc(W_x)$ 
and $\loc(W_y)$, so the points $x$ and $y$ can be connected using two curves in $W$ 
and a curve in $V^\sigma$. This implies that $X$ is rc$(\cW,V^\sigma)$-connected.\par

\medskip
Our next step consists in proving that $\rho_X=2$, showing that the numerical class 
of every irreducible component of any cycle in $\cW$ lies in the plane $\Pi$ spanned in $\cycl(X)$ by 
$[W]$ and $R_\sigma$. \\
Let $x \in X$ be a general point; by Lemma \ref{locy} we have 
$$\dim \loc(V^\sigma)_{\loc(W_x)} \ge l_\sigma +2i_X -2 \ge n-1,$$
therefore $E=\loc(V^\sigma)_{\loc(W_x)}$ and $\cycl(E)=\Pi$ by 
Lemma \ref{numequns}.\\
We have already proved that $-K_X \cdot W=4$ and $i_X=2$; therefore every reducible cycle in 
$\cW$ has exactly two irreducible components, and the families of deformations of
these components are unsplit. \\
Let $\Gamma_1 +\Gamma_2$ be a reducible cycle in $\cW$; without loss of generality
we can assume that $E \cdot \Gamma_1 >0$. Denote by $W^1$ a family of deformations of $\Gamma_1$;
being unsplit, the family $W^1$ cannot be dominating, hence 
for every $x \in \loc(W^1)$ we have $\dim \loc(W^1_x) \ge 2$ by Proposition \ref{iowifam}.
Since $E \cap \loc(W^1_x) \neq \emptyset$ it follows that $\dim (E \cap \loc(W^1_x)) \ge 1$ 
for every $x \in \loc(W^1)$, so $[W^1] \in \Pi$,
and consequently also $[W^2] \in \Pi$; it follows that $\rho_X=2$.\par

\medskip
Let now $T_Y$ be a minimal dominating family for $Y$ and let $T$ be the family
of deformations of the strict transform of a general curve in $T_Y$.
By Lemma \ref{doccia} we have $-K_X \cdot T = -K_{Y} \cdot T_Y \ge n-1$.\\
By this last inequality, the intersection $\loc(W_x) \cap \loc(T_x)$
for a general $x \in X$ has positive dimension; since $T$ is independent from 
$W$ -- recall that $E \cdot T=0$ and $E \cdot W >0$ --
the family $T$ cannot be locally quasi-unsplit.\\
Therefore, in the associated Chow family ${\mathcal T}$, there exists a reducible 
cycle $\Lambda=\Lambda_1 +\Lambda_2$
such that a family of deformations $T^1$ of $\Lambda_1$ is dominating and independent
from $T$.\\
The family $T^1$, being dominating, cannot be unsplit, hence $-K_X \cdot T^1 \ge 4$;
moreover, since $T^1$ is also independent from $T$ we have $E \cdot T^1>0$.
It follows that $E \cdot \Lambda_2 <0$ and so $-K_X \cdot \Lambda_2 \ge l_\sigma$
by Lemma \ref{Eneg2}. Therefore
$$-K_{Y} \cdot T_Y = -K_X \cdot T \ge l_\sigma + 2i_X = n+1$$
so $Y \simeq \proj^n$ by Corollary \ref{kedu}.\\
The center $B$ of $\sigma$ cannot be a linear subspace of $Y$, since otherwise
$i_X + l_\sigma =n+1$; take $l$ to be a proper bisecant of $B$ and  
let $\tl l$ be its strict transform: we have 
$$2=i_X \le -K_X \cdot \tl l = n+1 -2l_\sigma=  4 - l_\sigma,$$ 
hence $l_\sigma=2$ and $n=5$. \end{proof}

\begin{corollary} \label{nounsplit1}
Let $X$, $Y$, $R_\sigma$ and $E$ be as in \ref{setup2}.
Then there does not exist any family of rational curves $V$ independent from $R_\sigma$ such that 
$V_x$ is unsplit for some $x \in E$ and such that $E \subseteq \overline{\loc(V)}$.
\end{corollary}

\begin{proof}
Assume by contradiction that such a family exists.\\
First of all we prove that $V$ cannot be unsplit.
If this is the case, since on $X$ there are no unsplit dominating families it must be 
$\overline{\loc(V)}=\loc(V)=E$.\\
We can thus apply Lemma \ref{locy} a) to get that
$\dim \loc(V)_{F_\sigma}= n-1$ for every fiber $F_\sigma$ of $\sigma$. 
It follows that $E = \loc(V)_{F_\sigma}$ and therefore $\cone(E)=\twospan{R_\sigma}{[V]}$ by 
Lemma \ref{numequns}.\\
Since $V$ is a dominating unsplit family for the smooth variety $E$ we have 
$-K_E \cdot V = \dim \loc(V_x) + 1$, hence, by adjunction, 
$E \cdot V <0$; since $V$ is independent
from $R_\sigma$ it follows from Theorem \ref{Eneg1} 
that $Y$ is not a Fano manifold, a contradiction.\par

\medskip
Since $V$ is not unsplit we have $-K_X \cdot V \ge 2i_X$ and therefore,
for a point $x \in E$ such that $V_x$ is unsplit, we have  
$$\dim \loc(V_x) \ge -K_X \cdot V -1 \ge 2 i_X-1.$$ 
On the other hand, since $V$ is independent from $R_\sigma$, we have, 
for any fiber $F_\sigma$ of $\sigma$, that $\dim \loc(V_x) \cap F_\sigma \le 0$, 
hence $\dim \loc(V_x) \le n-l_\sigma = i_X +1$.\\
It follows that $i_X=2$, $-K_X \cdot V=4$ and $\dim \loc(V_x)=3$; the last
two equations, by Proposition \ref{iowifam} imply that 
$V$ is dominating.
Moreover, since $-K_X \cdot V=4$, the family $V$ is also locally unsplit, 
otherwise we would have a dominating family of lower degree, hence unsplit.\\
Since $E \cap \loc(V_x)$ is not empty and we cannot have $\loc(V_x) \subset E$ -- recall that
$V_x$ is unsplit and $V$ is independent from $R_\sigma$, so $\loc(V_x)$
can meet fibers of $\sigma$ only in points -- it follows that $E \cdot V>0$ and
we can apply Proposition \ref{p5}. \end{proof}

\begin{remark}\label{meetsB}
If $C_Y \subset Y$ is a curve which meets the center $B$ of the blow-up in $k$ points and is not 
contained in it, then $-K_{Y}\cdot C_Y \ge n-1 + (k-1)l_\sigma$.
\end{remark}

\begin{proof} \quad Let $C$ be the strict transform of $C_Y$: then the statement follows from the 
canonical bundle formula
$$-K_X = -\sigma^*K_{Y} - l_\sigma E,$$ 
which yields
$$-K_{Y} \cdot C_Y= -K_X \cdot C +l_\sigma E \cdot C 
\ge i_X +k l_\sigma  \ge n-1 + (k-1)l_\sigma .$$ 
\end{proof}

\begin{corollary}\label{redmeetsB} 
Let $W_Y$ be a minimal dominating family for $Y$
and assume that $-K_{Y} \cdot W_Y =n-1$. 
Assume that there exists a reducible cycle $\Gamma$ in 
$\cW_Y$ which meets $B$; then  $\Gamma \subset B$ and $\cone(B)=\onespan{[W_Y]}$.
\end{corollary}

\begin{proof} \quad Let $\Gamma_i$ be a component of $\Gamma$: we know that 
$-K_{Y} \cdot \Gamma_i < n-1$, so the whole cycle $\Gamma$ has to be contained in $B$
by remark \ref{meetsB}.\\
Let $W_Y^i$ be a family of deformations of $\Gamma_i$; the pointed locus $\loc(W_Y^i)_b$ is contained
in $B$ for every $b \in B$, again by remark \ref{meetsB}, hence
$$-K_{Y} \cdot W_Y^i \le \dim \loc(W_Y^i)_b \le \dim B = i_{X} \le i_{Y},$$
where the last inequality follows from \cite[Theorem 1, (iii)]{Bo}.\\
Therefore $T^i$ is unsplit and $B=\loc(W_Y^i)_b$, hence $\cone(B)=\onespan{[W_Y^i]}$.
It follows that all the components $\Gamma_i$ of $\Gamma$ are numerically proportional, and thus
they are all numerically proportional to $W_Y$.\end{proof}

We are now ready to prove the following

\begin{theorem}\label{noqu}
Let $X$ be a Fano manifold of dimension $n \ge 6$ and pseudoindex $i_X \ge 2$, which is the blow-up of 
another Fano manifold $Y$ along a smooth subvariety $B$ of dimension $i_X$;
assume that $X$ does not admit a quasi-unsplit dominating family of rational curves.
Then $Y \simeq \G(1,4)$ and $B$ is a plane of bidegree $(0,1)$. 
\end{theorem}

\begin{proof} The proof is quite long and complicated; we will divide it
into different steps, in order to make our procedure clearer.\par

\medskip
{\bf Step 1} \quad  \emph{A minimal dominating family of rational curves on $Y$ has 
anticanonical degree $n-1$.}\par

\medskip
Let $W_Y$ be a minimal dominating family of rational curves for
$Y$, and let $W$ be the family of deformations of the strict transform of 
a general curve in $W_Y$.\\
Apply \cite[Lemma 4.1]{AOlong} to $W$ (note that in the proof of that lemma
the minimality of $W$ is not needed). The first case in the lemma cannot occur by Corollary
\ref{nounsplit1}, so there exists a reducible 
cycle $\Gamma=\Gamma_\sigma + \Gamma_V +\Delta$ in $\cW$ with $[\Gamma_\sigma]$
belonging to $R_\sigma$, $\Gamma_V$ belonging to a family $V$, independent from $R_\sigma$, 
such that $V_x$ is unsplit for some $x \in E$, and $\Delta$ an effective rational 1-cycle. 
In particular 
\begin{equation}\label{W}
-K_X \cdot W \ge -K_X \cdot (\Gamma_\sigma + \Gamma_V +\Delta) \ge l_\sigma  +i_X \ge n-1.
\end{equation}
By the canonical bundle formula and Corollary \ref{E0} we have that
$$-K_{Y} \cdot W_Y = -K_X \cdot W  \ge n-1.$$
If  $-K_{Y} \cdot W_Y = n+1$ then $Y$ is a projective space by Corollary \ref{kedu}.
The center of $\sigma$ cannot be a linear subspace, otherwise $X$ would admit
an unsplit dominating family of rational curves; then if $l$ is a proper bisecant of $B$ and  
$\tl l$ is its strict transform we have 
$$2 \le i_X \le -K_X \cdot \tl l = n+1 -2l_\sigma=  4 - l_\sigma,$$ 
hence $l_\sigma=2$ and $n=5$, against the assumptions.\par

\smallskip
We can thus assume that $-K_{Y} \cdot W_Y \le n$.\\ 
Note that, by (\ref{W}), the reducible cycle $\Gamma$ 
has only two irreducible components $\Gamma_\sigma$ and $\Gamma_V$;
moreover the class of $\Gamma_\sigma$ is minimal in $R_\sigma$, hence 
$E \cdot \Gamma_\sigma=-1$, and $-K_X \cdot V \le i_X+1$. In particular $V$ is an unsplit
family.\\
Recalling that $E \cdot W = 0$ we get $E \cdot \Gamma_V = 1$.
Geometrically, a general curve in $V$ is the strict transform of a curve
in $W_Y$ which meets $B$ in one point; moreover, since a curve in $W_Y$ not contained in $B$
cannot meet $B$ in more than one point by Remark \ref{meetsB},
we have that 
\begin{equation}\label{lo}
\sigma(\loc(V) \setminus E)= \loc(W_Y)_B \setminus B.
\end{equation}
Assume that $-K_{Y} \cdot W_Y = n$; in this case $\rho_{Y}=1$ by Corollary \ref{minimaln}.\\
For a general point $y \in Y$, we have that $\loc(W_Y)_y$ is an effective, hence ample,
divisor, so it meets $B$. In particular we have $\dim \loc(W_Y)_B=n$, and
by (\ref{lo}) this implies that $V$ is dominating, against the assumptions
since $V$ is unsplit. This completes step 1.\par

\medskip
{\bf Step 2} \quad \emph{The strict transforms of curves in a minimal dominating family on $Y$ 
which meet $B$ fill up a divisor on $X$.}\par

\medskip
Let $x$ be a point in $E \cap \loc(V)$ and let $F_\sigma$ be the fiber of $\sigma$
containing $x$; since $\dim F_\sigma + \dim \loc(V_x) \le n$ we have
$$\dim \loc(V_x) \le n -l_\sigma  = i_X+1.$$
By inequality \ref{iowifam} we have that $\dim \loc(V) \ge n-2$; since $V$
is an unsplit family it cannot be dominating, so
we need to show that $\dim \loc(V) \not = n-2$.\\
Assume by contradiction that $\dim \loc(V)=n-2$; in this case, again by
inequality \ref{iowifam}, for every $x \in \loc(V)$ we have 
$\dim \loc(V_x)=i_X+1$, so for every $x \in X$ the intersection $\loc(V_x) \cap E$ dominates $B$.\\
Consider a point $x \in \loc(V) \setminus E$, denote by $y$ its image $\sigma(x)$ 
and consider $\loc(\cW_Y)_y$: since $\loc(V_x) \cap E$ dominates $B$, we have 
$B \subset \loc(\cW_Y)_y$. But cycles in $\cW_Y$ passing through 
$y$ and meeting $B$ are irreducible by corollary
\ref{redmeetsB}, so $B \subseteq \loc(W_Y)_y$ and by Lemma \ref{intmult} 
the numerical class of every curve 
in $B$ is an integral multiple of $[W_Y]$.
This fact together with Corollary \ref{redmeetsB} allows us to conclude that $B$ 
does not meet any reducible cycle in $\cW_Y$.\par

\medskip
We claim that a general curve $C$ of $W_Y$ is contained in the open subset $U$
of points $y \in Y$ such that $(W_Y)_y$ is proper. 
If this were not true, then $\loc(W_Y) \setminus U$ should have codimension one, and so there would
exist a family $W_Y^1$ of deformations of an irreducible component
of a cycle in $\cW_Y$ whose locus is a divisor; moreover this divisor
should have positive intersection number with $W_Y$.\\
This last condition would imply that $\loc(W_Y^1)$ has nonempty intersection
with $B$, since the numerical class of any curve in $B$ is
an integral multiple of $[W_Y]$, but we have proved that
$B$ does not meet any reducible cycle in $\cW_Y$, so we have reached a contradiction 
that proves the claim.\par

\medskip
Therefore we can apply Lemma \ref{minimaln-1} and get that $D:=\loc(W_Y)_{C}$ is a 
divisor and $\rho_{Y}=1$, since in the other case of the lemma we would find a family 
of anticanonical degree two meeting $B$, against Remark \ref{meetsB}.\\
Being $\rho_{Y}=1$ the effective divisor $D$ is ample, hence it meets $B$; 
therefore for a general curve $C$ in $W_Y$ there exists another curve in $W_Y$ which meets 
both $B$ and $C$; in other words, a general curve in $W_Y$ meets $\loc(W_Y)_B$, a contradiction 
since $\loc(W_Y)_B$ has codimension two in $Y$ by (\ref{lo}).\par

\medskip
{\bf Step 3} \quad \emph{The Picard number of $Y$ is one.}\par

\medskip
By (\ref{lo}) we have that $\dim \loc(W_Y)_B = \dim \loc(V)= n-1$.
This implies that $B$ contains curves whose numerical class is proportional to
$[W_Y]$, otherwise by Lemma \ref{locy} we would have  $\dim \loc(W_Y)_B= n$.\par

\smallskip
If  $B$ does not meet any reducible cycle in $\cW_Y$ we can argue as in the claim in step 2
and conclude that $\rho_{Y}=1$.\par

\smallskip
If else $B$ meets a reducible cycle in $\cW_Y$ then, by Corollary \ref{redmeetsB}, 
every curve in $B$ is numerically proportional to $[W_Y]$, hence
$\cone(\loc(W_Y)_B)=\langle [W_Y] \rangle$ and we conclude that $\rho_{Y}=1$ 
by Lemma \ref{minimalD}.\par

\medskip
{\bf Step 4} \quad \emph{The families of deformations of the strict transforms of curves in a minimal
dominating family on $Y$ which meet $B$ are extremal in $\cone(X)$.}\par

\medskip
Let $D =\loc(V)$; we have  $D \cdot W >0$, since 
$E \cdot W=0$ and $\pic(X)= \langle E, D \rangle$.
Therefore $\loc(W,V)_x=\loc(V)_{\loc(W_x)}$ is nonempty for a general $x \in X$,
and so has dimension $\ge n-2 + i_X -1 \ge n-1$ by Lemma \ref{locy}. It follows that
$i_X=2$ and $D=\loc(W,V)_x$.\\
The last equality, by Lemma \ref{numequns}, yields that  every curve in $D$ 
is numerically equivalent to a linear combination
$a[W] + b[V]$ with $a \ge 0$. \\
This implies that $\cone(D)$ is contained in the cone 
spanned by $[V]$ and by an extremal ray $R$ of $\cone(X)$. Since $E \cdot W =0$ and $E \cdot V > 0$
it must be $E \cdot R<0$, so $R = R_\sigma$ and $\cone(D)=\twospan{R_\sigma}{[V]}$.\\
Let $R_\tau$ be the extremal ray of $\cone(X)$ different from $R_\sigma$ and
denote by $\tau$ the associated contraction. The contraction $\tau$ is birational,
since $X$ does not admit quasi-unsplit dominating families of rational curves, therefore
its fibers have dimension at least two by inequality \ref{fiberlocus}. \\
We claim that $[V] \in R_\tau$; if we assume that this is not the case
then $D \cap \Exc(\tau)= \emptyset$, since $\cone(D)=\twospan{R_\sigma}{[V]}$.
In particular $D \cdot R_\tau=0$, so $D \cdot R_\sigma >0$ by \cite[Lemma 2.1]{BCW}.
By the same lemma the effective divisor $E$ is positive on $R_\tau$.\\
Let $F_\sigma$ and $F_\tau$ be two meeting fibers of the contractions $\sigma$ and $\tau$ respectively;
we have $\dim (F_\sigma \cap F_\tau) =0$, hence
$$n \ge \dim F_\sigma + \dim F_\tau \ge l_\sigma + l_\tau.$$
Therefore, recalling that $i_X=2$ and thus $l_\sigma=n-3$, we have $l_\tau \le 3$,
so $\dim \Exc(R_\tau)\ge n-2$ by inequality \ref{fiberlocus}.\\
In particular, if $F_\sigma$ is a fiber of $\sigma$ meeting $\Exc(\tau)$ we have
$$\dim (F_\sigma \cap \Exc(\tau)) \ge l_\sigma -2 \ge 1.$$
Let $C$ be a curve in $F_\sigma \cap \Exc(\tau)$; since $D \cdot R_\sigma >0$
we have $D \cap C \not = \emptyset$, hence $D \cap \Exc(\tau) \not = \emptyset$,
a contradiction that proves the extremality of $[V]$.\par

\medskip
{\bf Step 5} \quad \emph{The contraction of $X$ different from $\sigma$ is the blow-up
of $\proj^n$ along a smooth subvariety of codimension three.}\par

\medskip
Since $\Exc(\tau)= D=\loc(W,V)_x$ every fiber of $\tau$
is two-dimensional; we can apply \cite[Theorem 5.1]{AOsperays}
to get that $\tau\colon X \to Z$ is a smooth blow-up.\\
Let $T_Z$ be a minimal dominating family for $Z$ and $T^*$ a family of deformations
of the strict transform of a general curve in $T_Z$.\\
Among the families of deformations of the irreducible components of cycles in ${\mathcal T}^*$
there is at least one family which is dominating and locally unsplit; call it $T$.\\
Being dominating, $T$ cannot be quasi-unsplit, so we have $-K_X \cdot T \ge 4$,
hence for a general $x \in X$ we have $\dim \loc(T_x) \ge 3$.
Since $T$ is locally unsplit we also have $\cone(\loc(T_x))=\onespan{[T]}$.\\
On the other hand $\dim \loc(W_x) \ge n-2$,
so $\dim (\loc(T_x) \cap \loc(W_x)) \ge 1$. Therefore $T$ is numerically proportional
to $W$, since $\cone(\loc(W_x))=\onespan{[W]}$.\\ 
If $-K_X \cdot T < -K_X \cdot W$ then the images in $Y$ of the curves
in $T$ would be a dominating family for $Y$ of degree less than the degree of $W_Y$, a contradiction.
Therefore  $-K_{X} \cdot T \ge n-1$ and 
$$-K_{Z} \cdot T_Z = -K_X \cdot T^* \ge -K_X \cdot T + i_X \ge n+1,$$ 
so  $Z \simeq \proj^n$ by Corollary \ref{kedu} and $T_Z$
is the family of lines in $Z$.\par

\medskip
{\bf Step 6} \quad \emph{Conclusion.}\par

\medskip
Take $l_\sigma -2$ general sections  $H_i \in |\tau^*\Ol_{\proj^n}(1)|$; their intersection
${\mathcal I}$ is a Fano manifold of dimension five with two blow-up contractions of length two
$\sigma_{|{\mathcal I}} \colon {\mathcal I} \to Y'$ and 
$\tau_{|{\mathcal I}}\colon {\mathcal I} \to \proj^5$.\\
By the classification in \cite{CO} two cases are possible: either the center of 
$\tau_{|{\mathcal I}}$ is a Veronese 
surface or it is a cubic scroll contained in a hyperplane. The first case can be excluded
noting that, in our case, the degree of $E$ on a minimal curve in $R_\tau$ is one, since
$E \cdot W=0$ and $E \cdot R_\sigma =-1$.\\
It follows that $Y'$ is a del Pezzo manifold of degree five; $Y$ has $Y'$ as an ample section,
and therefore $Y$ is a del Pezzo manifold of degree five by repeated applications of 
\cite[Proposition A.1]{LPS}.
The only del Pezzo manifold of degree five and dimension greater than five is $\G(1,4)$.\end{proof}


\end{document}